\documentclass[a4paper,english,12pt]{amsart}

\usepackage{amsmath,amssymb}
\usepackage{amsthm}
\usepackage{enumerate}
\usepackage{babel}
\usepackage{a4wide}

\theoremstyle{plain}
\newtheorem{theorem}{Theorem}
\newtheorem*{theorem*}{Theorem}
\newtheorem{proposition}[theorem]{Proposition}
\newtheorem{corollary}[theorem]{Corollary}
\newtheorem*{corollary*}{Corollary}

\theoremstyle{remark}
\newtheorem{remark}[theorem]{Remark}
\theoremstyle{definition}

\DeclareMathSymbol{\Bbbk}{\mathalpha}{AMSb}{'174}
\newcommand{\BbbN}{\mathbb{N}}
\newcommand{\BbbZ}{\mathbb{Z}}

\newcommand{\BbbR}{\mathbb{R}}
\newcommand{\BbbC}{\mathbb{C}}
\newcommand{\domain}{\ensuremath{\Lambda}}
\newcommand{\field}{\ensuremath{\mathbf{k}}}
\newcommand{\graded}[2][]{G^{#1}(#2)}
\DeclareMathOperator{\grad}{gr}

\DeclareMathOperator{\GKdim}{GKdim}
\newcommand{\mono}[2]{\mathbf{#1}^{#2}}
\newcommand{\blackhole}[3]{#1^{-}_{#2}(#3)}
\newcommand{\mdeg}[2]{\mathrm{mdeg}_{#1}(#2)}
\newcommand{\esc}[2]{\langle #1, #2 \rangle}
\newcommand{\monodes}[3]{{#1}_{1}^{#2_1} \dots {#1}_{#3}^{#2_{#3}}}

\begin{document}

\title{Auslander-regular and Cohen-Macaulay Quantum groups}
\author{J. G\'omez-Torrecillas and F. J. Lobillo}
\address{Departamento de \'{A}lgebra, Universidad de Granada, E18071,
Granada, SPAIN} \email{torrecil@ugr.es} \email{jlobillo@ugr.es}
\maketitle

\section*{}

Let $U_q(C)$ be the quantum group or quantized enveloping algebra
in the sense of \cite{DeConcini/Kac:1990, DeConcini/Procesi:1993}
associated to a Cartan matrix $C$. A relevant property of $U_q(C)$
is that it can be endowed with a multi-filtration such that the
associated multi-graded algebra is an easy localization of the
coordinate ring of a quantum affine space \cite[Proposition
10.1]{DeConcini/Procesi:1993}. Thus, it is not surprising if we
claim that $U_q(C)$ is an Auslander-regular and Cohen-Macaulay
algebra (see, e.g., \cite{Bjork:1989} for these notions). However,
when one tries to construct a mathematically sound argument to
prove this, one realizes that there are not ready-to-use results
for this in the literature. Here we use re-filtering methods (see
Theorem \ref{eneauno}) similar to that in
\cite{Bueso/Gomez/Lobillo:1999unp} and
\cite{Bueso/Gomez/Lobillo:2000unp} to prove, in conjunction with
results from \cite{Bjork:1989} and \cite{McConnell/Stafford:1989},
that certain types of multi-filtered algebras are
Auslander-regular and Cohen-Macaulay (Theorem \ref{main}). This is
applied to obtain that $U_q(C)$ is Auslander-regular and
Cohen-Macaulay.

In this note, $K$ denotes a commutative ring and $\BbbN^n$ is the
free abelian monoid with $n$ generators $\epsilon_1, \dots,
\epsilon_n$. The elements in $\BbbN^n$ are vectors $\alpha =
(\alpha_1, \dots, \alpha_n)$ with non-negative integer entries. An
admissible order $\preceq$ on $\BbbN^n$ is a total order
compatible with the sum in $\BbbN^n$ and such that $0 \preceq
\alpha$ for every $\alpha \in \BbbN^n$. In this way, $\BbbN^n$
becomes a well-ordered monoid. A fundamental example of admissible
order on $\BbbN^n$ is the lexicographical order $\leq_{lex}$ with
$\epsilon_1 <_{lex} \cdots <_{lex} \epsilon_n$. Every vector
$\mathbf{w}$ with strictly positive entries gives an example of
admissible order $\preceq_{\mathbf{w}}$ by putting
\begin{equation}\label{ordenw}
\alpha \preceq_{\mathbf{w}} \beta \iff \begin{cases}
  \esc{\mathbf{w}}{\alpha} <  \esc{\mathbf{w}}{\beta} & \text{or} \\
 \esc{\mathbf{w}}{\alpha} =  \esc{\mathbf{w}}{\beta} \quad \text{and} \quad \alpha \leq_{\textrm{lex}}
 \beta
                        \end{cases}
\end{equation}
where $\esc{-}{-}$ denotes the usual dot product in $\BbbR^n$.

An $(\BbbN^n,\preceq)$--filtration on a $K$--algebra $R$ is a
family $F = \{ F_{\alpha}(R) ~|~ \alpha \in \BbbN^n \}$ of
$K$--submodules of $R$ such that
\begin{enumerate}
\item $F_{\alpha}(R) \subseteq F_{\beta}(R)$ for all
$\alpha \preceq \beta \in \BbbN^n$.
\item $F_{\alpha}(R) F_{\beta}(R) \subseteq F_{\alpha + \beta}(R)$
 for all $\alpha, \beta \in \BbbN^n$.
\item $\bigcup_{\alpha \in \BbbN^n} F_{\alpha}(R) = R$.
\item $1 \in F_{0}(R)$.
\end{enumerate}
The associated $\BbbN^n$--graded algebra is given by $G^F(R) =
\bigoplus_{\alpha \in \BbbN^n}G^F_{\alpha}(R)$, where
$G^F_{\alpha}(R) = F_{\alpha}(R)/\blackhole{F}{\alpha}{R}$ and
$\blackhole{F}{\alpha}{R} = \bigcup_{\beta \prec
\alpha}F_{\beta}(R)$. Further details can be found in
\cite{Gomez:1999}. The multi-degree of a nonzero element $r \in R$
is defined as $\mdeg{}{r} = \min \{ \alpha \in \BbbN^n ~|~ r \in
F_{\alpha}(R)\}$. When $n = 1$, the only admissible order is the
usual one and multi-filtrations are just positive filtrations. In
this case, the associated graded algebra will be denoted by
$\grad(R)$.

We will use extensively the following terminology: Let $\domain$
be a subalgebra of an algebra $R$, and let $x_1, \dots, x_n$ be
elements in $R$. A standard monomial in $x_1, \dots, x_n$ is an
expression $\mono{x}{\alpha} = \monodes{x}{\alpha}{n}$, where
$\alpha = (\alpha_1, \dots, \alpha_n) \in \BbbN^n$. Assume that an
element $r \in R$ can be written in the form
\begin{equation}\label{estandar}
r = \sum_{\alpha \in \BbbN^n} r_{\alpha}\mono{x}{\alpha} \qquad
(r_{\alpha }\in \domain)
\end{equation}
The expression \eqref{estandar} is called a (left) standard
representation of $r$. We will often refer as (left) polynomials
to the elements of $R$ having a standard representation.

\begin{theorem}\label{eneauno}
Let $\domain$ be a left noetherian subalgebra of a $K$--algebra
$R$, let $s$ be a positive integer and let $q_{ji} \in \domain$
for $1 \leqslant i < j \leqslant s$. The following statements are
equivalent
\begin{enumerate}[(i)]
\item \label{mfilt} There is an admissible order $\preceq$ on some $\BbbN^n$
and an $(\BbbN^n,\preceq)$--filtration $F = \{
  F_{\alpha}(R) ~|~ \alpha \in \BbbN^n \}$ on $R$ such that $F_0(R) = \domain$,
  every $F_{\alpha}(R)$ is finitely generated as a left
  $\domain$--module and $\graded[F]{R} = \domain [y_1; \sigma_1] \dots [y_s; \sigma_s]$
is an $\BbbN^n$--graded iterated Ore extension for some
homogeneous elements
  $y_1,\dots, y_s$ such that  $\sigma_j(y_i) =
  q_{ji}y_i$ for every $1 \leqslant i < j \leqslant s$.
\item \label{filt} There is an $\BbbN$--filtration $\{ R_n
  ~|~ n \in \BbbN \}$ on $R$ such that $R_0 = \domain$, every $R_n$ is finitely generated
  as a left $\domain$--module and $\grad(R) = \domain [y_1; \sigma_1] \dots [y_s; \sigma_s]$
  is an $\BbbN$--graded iterated Ore extension for some homogeneous elements  $y_1, \dots, y_s$
   such that $\sigma_j(y_i) = q_{ji} y_i$ for every $1 \leqslant i < j \leqslant s$.
\item \label{PBW} There are elements $x_1, \dots, x_s \in R$, an admissible order $\preceq'$ on $\BbbN^s$, and
finite subsets $\Gamma_{ji}, \Gamma_k \subseteq \BbbN^s$ for $1
\leqslant i < j \leqslant s, 1 \leqslant k \leqslant s$
  with $\max_{\preceq'}\Gamma_{ji} \prec' \epsilon_i +
  \epsilon_j$ and $\max_{\preceq'}\Gamma_{k} \prec' \epsilon_k$
  such that $\{
\mono{x}{\alpha} ~|~ \alpha \in \BbbN^s \}$ is a basis of $R$ as a
left $\domain$--module and
  $x_j x_i = q_{ji} x_i x_j + \sum_{\alpha \in \Gamma_{ji}} c_\alpha
  \mono{x}{\alpha}$ and for all $a \in \domain$, $x_k a = a^{(k)}x_k + \sum_{\alpha
\in \Gamma_i} c_\alpha \mono{x}{\alpha}$.
\end{enumerate}
\end{theorem}

\begin{proof}
\eqref{mfilt} implies \eqref{PBW}. Let $\alpha_i \in \BbbN^n$
denote the multi-degree of $y_i$ for $1 \leqslant i \leqslant s$.
Clearly, $\{\mono{y}{\gamma} ~|~ \gamma \in \BbbN^s \}$ is a basis
of $\graded[F]{R}$ as a left $\domain$--module. Thus, given $r \in
R$, the homogeneous element $r + \blackhole{F}{\mdeg{}{r}}{R} \in
\graded[F]{R}$ has a unique representation as homogeneous standard
left polynomial in $y_1, \dots, y_s$ with coefficients in
$\domain$. Thus,
\begin{equation}\label{homog}
r + \blackhole{F}{\mdeg{}{r}}{R} = \sum_{\gamma_1\alpha_1 + \cdots
+ \gamma_s\alpha_s = \mdeg{}{r}}c_{\gamma}\mono{y}{\gamma},
\end{equation}
where the $c_{\gamma}$'s are in $\domain$. Choose, for each $i =
1, \dots, s$, an element $x_i \in F_{\alpha_i}(R)$ such that $y _i
= x_i + \blackhole{F}{\alpha_i}{R}$. Let $M$ denote the $s \times
n$ matrix whose rows are $\alpha_1, \dots, \alpha_s$. Write the
equality \eqref{homog} as
\begin{equation}\label{neweq1}
r + \blackhole{F}{\mdeg{}{r}}{R} = \sum_{\gamma M =
\mdeg{}{r}}c_{\gamma}\mono{x}{\gamma} +
\blackhole{F}{\mdeg{}{r}}{R}
\end{equation}
Therefore, we can prove by induction on $\mdeg{}{r}$ that
\begin{equation}\label{neweq2}
r = \sum_{\gamma M \preceq \mdeg{}{r}}a_{\gamma}\mono{x}{\gamma},
\end{equation}
where $a_{\gamma} \in \domain$. To deduce that $\{
\mono{x}{\gamma} ~|~ \gamma \in \BbbN^s \}$ is a basis for
${}_\domain R$ we only need to check the linear independence.
Given a relation
\begin{equation}\label{ld}
\sum_{\gamma M \preceq \alpha} a_{\gamma}\mono{x}{\gamma} = 0,
\end{equation}
we proceed by induction on $\alpha$. The relation \eqref{ld} can
be written as
\begin{equation}\label{ld2}
\sum_{\gamma M = \alpha} a_{\gamma}\mono{x}{\gamma} + \sum_{\gamma
M \prec \alpha} a_{\gamma}\mono{x}{\gamma} = 0
\end{equation}
which, in $\graded[F]{R}$, gives
\begin{displaymath}
\sum_{\gamma M = \alpha} a_{\gamma}\mono{y}{\gamma} = 0
\end{displaymath}
As the monomials $\mono{y}{\gamma}$ are $\domain$--linearly
independent, we have that $a_{\gamma} = 0$ for $\gamma M =
\alpha$. The remaining coefficients are zero by induction in view
of \eqref{ld2}.

Let $a \in \domain$ and $i \in \{1, \dots, s \}$. Since
$G_0^{F}(R) = F_0(R) = \domain$ and $y_ia = \sigma_i(a)y_i$ we get
$\sigma_i(a)$ has degree $0$, i.e., $\sigma_i(a) \in \domain$.
Write $a^{(i)} = \sigma_i(a)$. Then
\begin{equation}\label{conmui}
0 = y_ia - a^{(i)}y_i = (x_ia - a^{(i)}x_i) +
\blackhole{F}{\alpha_i}{R}
\end{equation}
Since $\domain$ is left noetherian and $F_{\alpha_i}(R)$ is
finitely generated as a left $\domain$--module, we have that
$\blackhole{F}{\alpha_i}{R}$ is a noetherian left
$\domain$--module. Thus, we deduce from \eqref{conmui}, in
conjunction with \eqref{neweq2}, that
\begin{equation}\label{neweq3}
x_i a = a^{(i)}x_i + \sum_{\gamma  \in \Gamma_i}a_{\gamma}\mono{x}{\gamma},
\end{equation}
for some $a_{\gamma} \in \domain$, where $\Gamma_i$ is a finite
subset of $\BbbN^s$ such that $\gamma M \prec \alpha_i$ for every
$\gamma \in \Gamma_i$. On the other hand, for $1 \leq i < j \leq
s$, we have
\begin{multline*}
0 = y_jy_i - q´_{ji}y_iy_j \\ = (x_j +
\blackhole{F}{\alpha_j}{R})(x_i + \blackhole{F}{\alpha_i}{R}) -
q_{ji}(x_i + \blackhole{F}{\alpha_i}{R})(x_j + \blackhole{F}{\alpha_j}{R}) \\
= (x_jx_i - q_{ji}x_ix_j) + \blackhole{F}{\alpha_i + \alpha_j}{R},
\end{multline*}
which entails, by \eqref{neweq2},
\begin{equation}\label{neweq4}
x_jx_i - q_{ji}x_ix_j = \sum_{\gamma  \in
\Gamma_{ij}}a_{\gamma}\mono{x}{\gamma},
\end{equation}
where $\Gamma_{ij}$ is a finite subset of $\BbbN^s$ such that
$\gamma M \prec \alpha_i + \alpha_j$ for every $\gamma \in
\Gamma_{ij}$. Let $\preceq'$ be the admissible order on $\BbbN^s$
defined by
\begin{equation}\label{ordenM}
\gamma \preceq' \mu \iff \begin{cases}
  \gamma M \prec \mu M & \text{or} \\
  \gamma M = \mu M \quad \text{and} \quad \gamma \leq_{\textrm{lex}} \mu
                        \end{cases}
\end{equation}
Since $\alpha_i = \epsilon_iM$ for every $i = 1, \dots, s$, the
relations \eqref{neweq3} and \eqref{neweq4} can be written as
\begin{equation}\label{cadai}
x_ia-a^{(i)}x_i = \sum_{\gamma \prec' \epsilon_i \atop \gamma \in
\Gamma_i }a_{\gamma} \mono{x}{\gamma}
\end{equation}
and
\begin{equation}\label{cadaij}
x_jx_i - q_{ji}x_ix_j = \sum_{\gamma \prec' \epsilon_i +
\epsilon_j \atop \gamma \in \Gamma_{ij}}a_{\gamma}\mono{x}{\gamma},
\end{equation}
which gives \eqref{PBW}.

\eqref{PBW} implies \eqref{filt}. First, notice that, by
hypothesis, the relations \eqref{cadai} and \eqref{cadaij} are
satisfied. Let
\begin{displaymath}
C = \{0\} \cup \Big( \bigcup_{1 \leq i \leq s} C_i \Big) \cup \Big( \bigcup_{1
  \leq i < j \leq s} C_{ij} \Big),
\end{displaymath}
where $C_i = \Gamma_i - \epsilon_i$ and $C_{ij} = \Gamma_{ij} -
\epsilon_i - \epsilon_j$. Clearly $C$ is a finite subset of
$\BbbZ^s$ whose maximum with respect to $\preceq$ is $0$. By
\cite[Corollary 2.2]{Bueso/Gomez/Lobillo:1999unp} (see also
\cite{Mora/Robbiano:1988} and \cite{Weispfenning:1987}), there is
$\mathbf{w} = (w_1, \dots, w_s) \in \BbbN^s_+$ such that
$\esc{\mathbf{w}}{\alpha} < 0$ for every $\alpha \in C$. This
implies that the relations \eqref{cadai} and \eqref{cadaij} can be
written as
\begin{equation}\label{cadaiw}
x_ia - a^{(i)}x_i = \sum_{\esc{\mathbf{w}}{\gamma} <
  w_i}a_{\gamma}\mono{x}{\gamma}
\end{equation}
and
\begin{equation}\label{cadaijw}
x_jx_i - q_{ji}x_ix_j = \sum_{\esc{\mathbf{w}}{\gamma} < w_i +
w_j}a_{\gamma}\mono{x}{\gamma}
\end{equation}
By \cite[Proposition 1.13]{Bueso/Gomez/Lobillo:1999unp}, the
family $\{H_{\alpha}(R) ~|~ \alpha \in \BbbN^s\}$ where
$H_\alpha(R)$ is the left \domain--module generated by the set
$\{\mono{x}{\beta} ~|~ \beta \preceq_\mathbf{w} \alpha \}$, is an
$(\BbbN^s,\preceq_\mathbf{w})$--filtration on $R$. Since
$\mathbf{w}$ has no zero component, it follows that
$H_{\alpha}(R)$ is finitely generated as a left $\domain$--module
for every $\alpha$. For each $n \in \BbbN$, define $R_n =
\bigcup_{\esc{\mathbf{w}}{\alpha} \leqslant n}H_{\alpha}(R)$,
which is a finitely generated left $\domain$--module. A
straightforward verification shows that $\{R_n ~|~ n \in \BbbN\}$
is a filtration on $R$. Clearly, $R_n =
\sum_{\esc{\mathbf{w}}{\alpha} \leqslant n}\domain
\mono{x}{\alpha}$ for every $\alpha \in \BbbN^s$. Finally, let
$y_i = x_i + R_{w_i - 1}$ for $1 \leqslant i < j \leqslant s$. By
\eqref{cadaiw} and \eqref{cadaijw}, $x_ia = a^{(i)}x_i$ for every
$a \in \domain$ and $x_jx_i = q_{ji}x_ix_j$. Moreover, since the
monomials $\mono{x}{\alpha}$ are $\domain$--linearly independent,
it follows that $\{\mono{y}{\alpha} ~|~ \alpha \in \BbbN^s\}$ is a
left \domain--basis for $\grad(R)$. It follows from
\cite[2.1.(iii)]{Goodearl/Letzter:1994} that
\begin{displaymath}
\grad(R) \cong \domain[y_1; \sigma_1]\cdots[y_s;\sigma_s]
\end{displaymath}
Finally, \eqref{filt} implies \eqref{mfilt} obviously.
\end{proof}

In the following corollary, $K_0(R)$ denotes the Grothendieck
group of $R$. Of course, the corollary says something new for
rings satisfying \eqref{mfilt} or \eqref{PBW} in Theorem
\ref{eneauno}.

\begin{corollary}\label{regular}
Assume $R$ satisfies one (and then all) of the equivalent
conditions of Theorem \ref{eneauno}. Suppose, in addition, that
$\domain$ is right noetherian, $q_{ji}$ is a unit of $\domain$ for
$1 \leq i < j \leq s$ and that $\sigma_i$ is an automorphism of
 $\domain$ for $i=1, \dots, s$.
\begin{enumerate}
\item
If every cyclic right $\domain$--module has finite projective
dimension, then $K_0(\domain) \cong K_0(R)$.
\item
If $\Lambda$ is Auslander-regular then $R$ Auslander-regular.
\end{enumerate}
\end{corollary}
\begin{proof}
The first statement is a consequence of \cite[Theorem
12.6.13]{McConnell/Robson:1988}. If $\Lambda$ is
Auslander-regular, then, by \cite[Theorem 4.2]{Ekstrom:1989},
$\grad(R) = \domain[y_1;\sigma_1] \cdots [y_s;\sigma_s]$ is
Auslander-regular. The result follows now from \cite[Theorem
3.9]{Bjork:1989}.
\end{proof}

\begin{theorem}\label{main}
Assume that $R$ is an algebra over a field $\field$ satisfying one
(and then all) of the equivalent conditions of Theorem
\ref{eneauno}. Suppose, in addition, that
\begin{enumerate}[(a)]
\item The scalars $q_{ji}$ are units of $\field$ and the
endomorphisms $\sigma_i : \domain \rightarrow \domain$ are
automorphisms.
\item $\domain$ is generated as an algebra by elements $z_1, \dots, z_t$
such that the standard filtration $\domain_n$ obtained by giving
degree $1$ to each $z_i$ satisfies that $gr(\domain)= \oplus_{n
\geqslant 0}\domain_n/\domain_{n-1}$ is a finitely presented and
noetherian algebra over $\field$.
\item $\sigma_i(\domain_1) \subseteq \domain_1$, for $i = 1, \dots, s$.
\item either $gr(\domain)$ or
$\domain[y_1;\sigma_1]\cdots[y_s;\sigma_s]$ is an
Auslander-regular and Cohen-Macaulay algebra.
\end{enumerate}
Then $R$ is an Auslander-regular and Cohen-Macaulay algebra.
\end{theorem}

\begin{proof}
Let $R_n$ be the filtration on $R$ given by Theorem \ref{eneauno}
with $\grad(R) = \domain[y_1;\sigma_1]\cdots[y_s;\sigma_s]$. Since
$\sigma_i(\domain_1) \subseteq \domain_1$ for every $i = 1, \dots,
s$ and the filtration $\domain_n$ is standard, we get that
$y_i\domain_n \subseteq \domain_ny_i$ for every $i = 1, \dots, s$
and every $n \geqslant 0$. Therefore, $\domain \subseteq
\domain[y_1;\sigma_1] \cdots [y_s;\sigma_s]$ is a
$\preceq_\mathbf{w}$--bounded extension of $\domain$ in the sense
of \cite[Definition 1.8]{Bueso/Gomez/Lobillo:1999unp}. Here,
$\mathbf{w} = (w_1, \dots, w_s)$ with $w_i = \deg(y_i)$, $i= 1,
\dots, s$. Let $\preceq$ be the admissible order defined by
\[
(i,\alpha) \preceq (j,\beta) \iff \begin{cases}
\alpha \prec_{\mathbf{w}} \beta & \text{or} \\
\alpha = \beta \text{ and } i \leq j &
\end{cases}
\]
Write $H(i,\alpha) = \sum_{(j,\beta) \preceq (i,\alpha)} \domain_j
\mono{y}{\beta}$. By \cite[Proposition
1.13]{Bueso/Gomez/Lobillo:1999unp}, these vector subspaces form a
$(\BbbN^{s+1},\preceq)$--filtration for $\grad{R}$. Let
$\grad(R)_{(n)} = \sum_{i+\esc{\mathbf{w}}{\alpha}\leq n}
\domain_i \mono{y}{\alpha}$. Since $$\grad(R)_{(n)} =
\bigcup_{\esc{(1,\mathbf{w})}{(i,\alpha)}}H_{(i,\alpha)},$$ it
follows that $\{\grad(R)_{(n)} ~|~ n \in \BbbN\}$ is a filtration
on $\grad(R)$. Moreover, the inclusion $\domain \subseteq
\grad(R)$ is a strict filtered morphism, hence $\grad(\domain)$
can be viewed as a subalgebra of $\grad(\grad(R))$. Therefore,
$\grad(\grad(R)) \cong
\grad(\domain)[y_1;\sigma_1]\cdots[y_s;\sigma_s]$. Here,
$\sigma_i$ denotes the graded automorphism induced by the
homonymous filtered automorphism of $\domain [y_1; \sigma_1]
\cdots [y_{i-1};\sigma_{i-1}]$. Since $\grad(\domain)$ is a
finitely presented and noetherian algebra, we see that
$\grad(\grad(R))$ enjoys the same properties. Thus, the filtration
$R_n$ satisfies the hypotheses  of \cite[Theorem
1.3]{McConnell/Stafford:1989}. Now every finitely generated left
$R$--module is endowed with a filtration such that $\grad(M)$ is
finitely generated. By \cite[Theorem
1.3]{McConnell/Stafford:1989}, $\GKdim (M) = \GKdim (\grad(M))$.
In particular, $\GKdim (R) = \GKdim (\grad(R))$. On the other
hand, from the proof of \cite[Theorem 3.9]{Bjork:1989} we obtain
that $j_R(M) = j_{\grad(R)}(\grad(M))$. If we assume that
$\grad(R) = \domain [y_1;\sigma_1]\cdots[y_s;\sigma_s]$ is
Cohen-Macaulay, then
\begin{displaymath}
\mathrm{GKdim}(R) = \mathrm{GKdim} (\grad(R)) = j_{\grad(R)}(\grad(M)) +
\GKdim (\grad(M)) = j_R(M) + \GKdim (M),
\end{displaymath}
whence $R$ is Cohen-Macaulay too.

Lastly, if $\grad(\domain)$ is Cohen-Macaulay, then
$\grad(\grad(R))$ satisfies the hypotheses of
\cite[Lemma]{Levasseur/Stafford:1993}, which implies that it is
Cohen-Macaulay. Since the filtration $\grad(R)_{(n)}$ is
finite-dimensional, we obtain that $\grad(R)$ is Cohen-Macaulay.
Thus, $R$ is Cohen-Macaulay by the foregoing argument.
\end{proof}

If $Q = (q_{ij})$ is a multiplicatively anti-symmetric $s \times
s$ matrix with coefficients in $\field$, the coordinate ring of
the quantum affine space $\mathcal{O}_{Q}(\field^s) =
\field_Q[x_1, \dots, x_s]$ is the $\field$--algebra generated by
$x_1, \dots, x_s$ subject to the relations $x_jx_i =
q_{ji}x_ix_j$.

For our purposes, we are interested in certain localizations of
$\mathcal{O}_{Q}(\field^s)$. Thus, consider some of the variables
which, for simplicity, we assume to be $x_1, \dots x_t$ with $t
\leqslant s$. Since $x_1, \dots, x_t$ are normal elements, they
generate a multiplicatively closed Ore set, so that we can
construct the localized algebra
\[
\mathbf{k}_{Q}[x_1^{\pm 1}, \dots, x_t^{\pm 1}, x_{t+1}, \dots,
x_{s}]
\]

Although the following proposition should be well-known, we have
not found a precise reference.

\begin{proposition}\label{grCM}
The algebra $A = \mathbf{k}_{Q}[x_1^{\pm 1}, \dots, x_t^{\pm 1},
x_{t+1}, \dots, x_{s}]$ is Auslander-regular and Cohen-Macaulay.
\end{proposition}

\begin{proof}
Clearly, $A$ is an iterated Ore extension of a McConnell-Pettit
algebra, whence its global homological dimension is finite by
\cite[3.1]{McConnell/Pettit:1988} and \cite[Theorem
4.2]{Ekstrom:1989}. On the other hand,
\[
\mathbf{k}_{Q}[x_1, \dots, x_t,x_{t+1}, \dots, x_s]
\]
is Auslander-regular and Cohen-Macaulay (see, e.g., \cite[Theorem
3.5]{Goodearl/Lenagan:1996}). By \cite[Proposition
2.1]{Ajitabh/Smith/Zhang:1999}, $A$ satisfies the Auslander
condition. Since the multiplicative set generated by $x_1, \dots,
x_t$ consists of monomials, which are local normal elements, we
have, by \cite[Theorem 2.4]{Ajitabh/Smith/Zhang:1999}, that our
algebra $A$ is Cohen-Macaulay.
\end{proof}

\begin{theorem}
The quantized enveloping $\BbbC (q)$--algebra $U_q(C)$ associated
to a Cartan matrix $C$ is Auslander-regular and Cohen-Macaulay.
\end{theorem}

\begin{proof}
Accordingly with \cite[Proposition 10.1]{DeConcini/Procesi:1993},
$U = U_q(C)$ is endowed with a
$(\mathbb{N}^{n},\preceq)$--filtration $\{ F_{\alpha}(U) ~|~
\alpha \in \BbbN^n \}$ for some $n$ and a lexicographical order
$\preceq$ in such a way that the multi-graded associated algebra
$\graded[F]{U} \cong \mathbb{C}(q)_{Q}[x_1^{\pm 1}, \dots,
x_t^{\pm 1}, x_{t+1}, \dots, x_{s}]$ for a certain
multiplicatively anti-symmetric matrix $Q$. By Proposition
\ref{grCM}, $\graded[F]{U}$ is Auslander-regular and
Cohen-Macaulay. Moreover, $F_0(U) = \mathbb{C}(q)[z_1^{\pm
1},\dots,z_t^{\pm 1}]$, a commutative Laurent polynomial ring.
Filter $F_0(U)$ with the standard filtration obtained by giving
degree $1$ to $z_i^{\pm1}$ ($i = 1, \dots, t$). Then $gr(F_0(U))$
is a factor algebra of the commutative polynomial ring in $2t$
variables with coefficients in $\mathbb{C}(q)$. In particular, it
is finitely presented and noetherian. Therefore, the hypotheses of
Theorem \ref{main} are fulfilled and, hence, $U_q(C)$ is
Auslander-regular and Cohen-Macaulay.
\end{proof}

\begin{remark}
In \cite[Proposition 2.2]{Brown/Goodearl:1997} it is shown that
$U_q(C)$ is Auslander-regular. It is also proved \cite[Theorem
2.3]{Brown/Goodearl:1997} that $U_q(C)$ is Cohen-Macaulay with
respect to the Krull dimension in case $q$ is a root of unity..
\end{remark}

\begin{remark}
The normal separation of the prime spectrum of $U_q(C)$ would
imply in view of Theorem \ref{grCM} and \cite[Theorem
1.6]{Goodearl/Lenagan:1996} that $U_q(C)$ is catenary. However,
the (classical) universal enveloping algebras are not normally
separated in general. So, as the referee pointed out, it is
interesting to know if $U_q(C)$ does not really enjoy this
property and why.
\end{remark}

\providecommand{\bysame}{\leavevmode\hbox to3em{\hrulefill}\thinspace}

\end{document}